\documentclass{article} 

\usepackage{amsthm}
\newtheorem{proposition}{Proposition}[section]
\newtheorem{theorem}{Theorem}[section]
\newtheorem{remark}{Remark}[section]

\usepackage{makeidx}         
\usepackage{graphicx}        
\usepackage{multicol}        
\usepackage[bottom]{footmisc}

\usepackage{newtxtext}       %
\usepackage[varvw]{newtxmath}       

\usepackage{hyperref}
\usepackage{cprotect}

\makeindex             

\begin{document}

\title{A New Diffusive Representation for Fractional Derivatives and its Application}

\author{Kai Diethelm
\footnote{Faculty of Applied Natural Sciences and Humanities, University of Applied Sciences Würzburg-Schweinfurt, 
	Ignaz-Schön-Str.\ 11, 97421 Schweinfurt, Germany, email{kai.diethelm@fhws.de}}}
\maketitle

\abstract{Diffusive representations of fractional derivatives have proven to be useful tools 
	in the construction of fast and memory efficient numerical methods for solving fractional
	differential equations. A common challenge in many of the known variants of this approach is that
	they require the numerical approximation of some integrals over an unbounded integral whose
	integrand decays rather slowly which implies that their numerical handling is difficult and costly.
	We present a novel variant of such a diffusive representation. This form also requires the numerical
	approximation of an integral over an unbounded domain, but the integrand decays much faster.
	This allows to use well established quadrature rules with much better convergence properties.}

\section{Introduction and Statement of the Problem}

\subsection{Classical Discretizations in Fractional Calculus}

The efficient numerical solution of initial value problems with fractional differential equations 
like, e.g., 
\begin{equation}
	\label{eq:ivp}
	D^\alpha_a y(t) = f(t, y(t)), \qquad y(a) = y_0,
\end{equation}
is a significant 
computational challenge due to, among other reasons, the non-locality of fractional
differential operators. 
In our formulation \eqref{eq:ivp}, $D^\alpha_a$ denotes the standard
Caputo differential operator of order $\alpha$ with starting point 
$a \in \mathbb R$ \cite[Chapter 3]{Di2010},
and we assume here and throughout some other parts of this paper that $0 < \alpha < 1$
(although we explicitly point out that the generalization of our findings to the case
that $\alpha$ is a noninteger number greater than $1$ is a relatively straightforward matter).

When dealing with the problem \eqref{eq:ivp}, one usually introduces
a discretization of the interval $[a, a+T]$, say, on which the solution is sought by
defining some grid points $a = t_0 < t_1 < t_2 < \cdots < t_N = a+T$. For each
grid point $t_j$, $j = 1, 2, \ldots, N$, typical numerical methods then introduce an 
approximation formula for a discretization of $D_a^\alpha y(t_j)$ based on function values of
$y$ at the grid points, replace the exact fractional derivative in eq.~\eqref{eq:ivp} by this
approximation, discard the approximation error and solve the resulting algebraic equation to
obtain an approximation for $y(t_j)$. In their standard forms, classical methods like
fractional linear multistep methods \cite{Lub:lmm-abel,Lub:discr-fc} or the Adams method
\cite{DFF2002,DFF2004} require $O(j)$ operations to compute the required approximation at the
$j$-th grid point, thus leading to an $O(N^2)$ complexity for the overall calculation of the 
approximate solution at all $N$ grid points. Moreover, the construction of the algorithms
requires the entire history of the process to be in the active memory at any time, thus leading
to an $O(N)$ memory requirement. This may be prohibitive in situations like,
e.g., the simulation of the mechanical behaviour of viscoelastic materials via some 
finite element code where a very large number of such differential equations needs 
to be solved simultaneously \cite{HSL2019}.

Numerous modifications of these basic algorithms have been proposed to resolve these
issues. Specifically (see, e.g., \cite[Section 3]{DKLMT}), one may use FFT techniques 
to evaluate the sums that arise in the formulas \cite{Ga2018,HLS1985,HLS1988}, thus reducing
the overall computational complexity to $O(N \log^2 N)$; however, this approach does not improve the
memory requirements. Alternatively, nested mesh techniques \cite{DF2006,FS} can be employed;
this typically reduces the computational complexity to $O(N \log N)$, and some of these methods are
also able to cut down the active memory requirements to $O(\log N)$.

\subsection{Diffusive Representations in Discretized Fractional Calculus} 

From the properties recalled above, it becomes clear that none of the
schemes mentioned so far allows to reach the level known for traditional algorithms for 
first order initial value problems that, due to their differential operators being local, 
have an $O(N)$ complexity and an $O(1)$ memory requirement. However, it is possible 
to achieve these perfomance features by using methods based on diffusive representations
for the fractional derivatives \cite{Mo1998}. Typically, such representations take the form
\begin{equation}
	\label{eq:diff-rep}
	D_a^\alpha y(t) = \int_0^\infty \phi(w, t) \mathrm d w
\end{equation}
where, for a fixed value of $w$, the function $\phi(w, \cdot)$ is characterized as the solution
to an initial value problem for a first order differential equation the formulation of which contains
the function $y$ whose fractional derivative is to be computed. In the presently available literature,
many different special cases of this representation are known, e.g.\ the version of Yuan and
Agrawal \cite{YA} (originally proposed in that paper for $0 < \alpha < 1$ and extended to
$1 < \alpha < 2$ in \cite{TR} and to general positive noninteger values of $\alpha$ in
\cite{Di2008}; see also \cite{SG} for further properties of this method) where the associated 
initial value problem reads
\begin{subequations}
	\label{eq:ya}
	\begin{equation}
		\label{eq:ya-ode}
		\frac{\partial \phi^{\mathrm{YA}}}{\partial t} (w,t) \!
		        = \!- w^2 \phi^{\mathrm{YA}}\!(w,t) 
			         \! + \! (\!-1\!)^{\lfloor \!\alpha\! \rfloor} \frac{2 \sin \pi \alpha} \pi 
			                w^{2\alpha-2\lceil \! \alpha \! \rceil+1}  y^{(\lceil\! \alpha\! \rceil)}(t), 
		\,\,
		\phi^{\mathrm{YA}}\!(w, a) = 0,
	\end{equation}
  	such that the function $\phi^{\mathrm{YA}}$ has the form
	\begin{equation}
		\label{eq:ya-phi}
		\phi^{\mathrm{YA}}(w,t) 
			= (-1)^{\lfloor \alpha \rfloor} \frac{2 \sin \pi \alpha} \pi 
					w^{2\alpha-2\lceil \alpha \rceil+1} \int_a^t y^{(\lceil \alpha \rceil)}(\tau) 
				              \exp(-(x-\tau) w^2) \mathrm d \tau.
	\end{equation}
\end{subequations}
An alternative has been proposed by Chatterjee \cite{Chatt} (see also \cite{SC}) using the
initial value problem
\begin{subequations}
	\label{eq:sc}
	\begin{equation}
		\label{eq:sc-ode}
		\frac{\partial \phi^{\mathrm{C}}}{\partial t} (w,t) \!
		        = \! - w^{1 / (\alpha - \lceil \! \alpha \! \rceil + 1)} \phi^{\mathrm{C}}\!(w,t) 
			          + (\!-1\!)^{\lfloor \! \alpha \! \rfloor} \frac{\sin \pi \alpha} {\pi  (\alpha - \lceil \!\alpha\!\rceil + 1)}
			                y^{(\lceil \!\alpha\! \rceil)}(t), 
		\quad
		\phi^{\mathrm{C}}\!(w, a) = 0,
	\end{equation}
  	such that the function $\phi^{\mathrm{C}}$ has the form
	\begin{equation}
		\label{eq:sc-phi}
		\phi^{\mathrm{C}}(w,t) 
    			=  \frac{(-1)^{\lfloor  \alpha  \rfloor}\sin \pi  \alpha}{\pi ( \alpha  - \lceil  \alpha  \rceil + 1)} 
				      \int_a^t y^{(\lceil  \alpha  \rceil)}(\tau) 
					     \exp\left(-(t-\tau) w^{1/( \alpha  - \lceil  \alpha  \rceil + 1)}\right) 
   					  \mathrm d \tau.
   	\end{equation}
\end{subequations}

In either case (or in the case of the many variants thereof that have been proposed; 
cf., e.g., \cite{Ba2019,BS,Li,McL,ZCSHN}), the numerical calculation of $D^\alpha_a y(t_j)$ requires 
\begin{enumerate}
\item a quadrature formula 
	\begin{equation}
		\label{eq:qf1}
		\sum_{k=1}^K \lambda_k \phi(w_k, t_j) \approx  \int_0^\infty \phi(w, t_j) \mathrm d w
			= D_a^\alpha y(t_j) 
	\end{equation}
	with nodes $w_1, w_2, \ldots, w_K \in [0, \infty)$ and weights $\lambda_1, \lambda_2, 
	\ldots, \lambda_K \in \mathbb R$	for numerically evaluating the integral in eq.~\eqref{eq:diff-rep}
\item a standard numerical solver for the associated differential equation (e.g., a linear multistep
	method) to approximately compute, for each $k \in \{ 1, 2, \ldots, K \}$, the values $\phi(w_k, t_j)$ 
	required to evaluate the formula \eqref{eq:qf1}.
\end{enumerate}
Evidently, the run time and the memory requirements of the operation in step~1 do not depend on $j$.
Also, one can perform step 2 in an amount of time that is independent of $j$. Furthermore, if an $\ell$-step 
method is used in step 2, one needs to have (approximate) information about $y(t_{j-1}), y(t_{j-2}), 
\ldots y(t_{j-\ell})$ which has to be kept in the active memory---but the amount of storage space required 
for this purpose is also independent of $j$.

In summary, approaches of this type require $O(1)$ arithmetic operations per time step, i.e.\ we
have a computational cost of $O(N)$ for all $N$ time steps combined, and the required amount of 
memory is $O(1)$ as desired. A further pleasant property of these methods is that they
impose no restrictions at all on the choice of the grid points $t_j$ whereas this can not always be
achieved with the other approaches. 
Thus, from a theoretical point of view, algorithms of this form 
are very attractive. In practice, however, the implied constants in the $O$-terms may be very large.
This is due to the following observation \cite[Theorems 3.20(b) and 3.21(b)]{Di2010}:
\begin{proposition}
	\label{prop:asymp-yac}
	Let $t \in [a, a+T]$ be fixed. Then, for $w \to \infty$, we have
	\[
		\phi^{\mathrm{YA}}(w, t) = c^{\mathrm{YA}} w^{q_{\mathrm{YA}}} (1 + o(1))
		\quad \mbox{ with } \quad
		q_{\mathrm{YA}} = 2 \alpha - 2 \lceil\alpha\rceil - 1 \in (-3, -1)
	\]
	and
	\[
		\phi^{\mathrm{C}}(w, t) = c^{\mathrm{C}} w^{q_{\mathrm C}} (1 + o(1))
		\quad \mbox{ with } \quad
		q_{\mathrm{C}} = - \frac 1 {\alpha - \lceil\alpha\rceil + 1} < -1
	\]
	where $c^{\mathrm{C}}$ and $c^{\mathrm{YA}}$ are some constants independent of $w$ 
	(that may, however, depend on $t$, $a$, $\alpha$ and $y$).
\end{proposition}
\vskip-0.5cm
\begin{figure}[htb]
	\centering
	\includegraphics[width=0.7\textwidth]{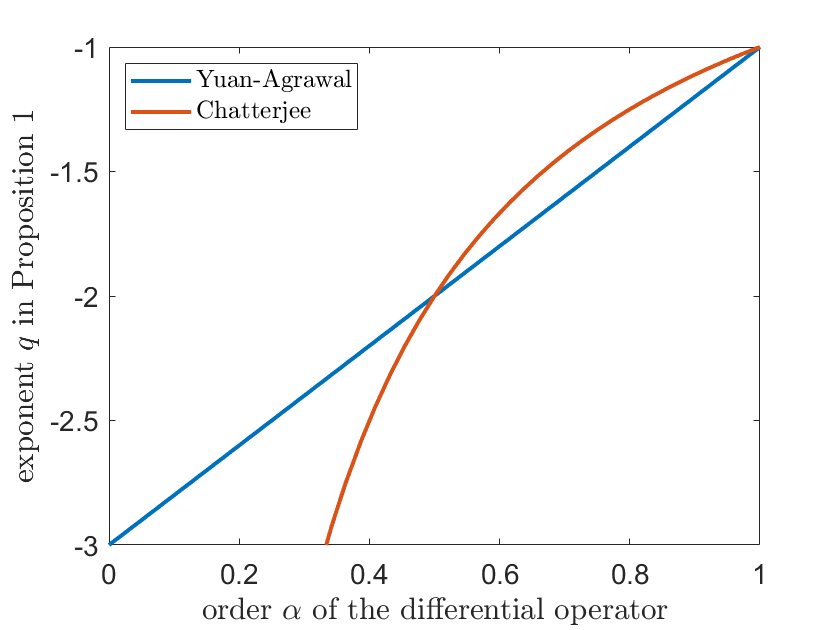}
	\caption{\label{fig:q-vs-alpha}Behaviour of the exponents 
		$q_{\mathrm{YA}}$ (blue) and $q_{\mathrm{C}}$ (orange) 
		as introduced in Proposition \ref{prop:asymp-yac} vs.\ the 
		order $\alpha$ of the differential operator}
\end{figure}

From Proposition \ref{prop:asymp-yac} we can see that the integrands in eq.\ \eqref{eq:diff-rep}
decay to zero in an algebraic manner as $w \to \infty$.
Figure \ref{fig:q-vs-alpha} shows the behaviour of the exponents $q_{\mathrm{YA}}$ and 
$q_{\mathrm C}$ as they depend on $\alpha$. It can be seen that the exponents are less than $-1$ for
all $\alpha \in (0,1)$. This suffices to assert that the integrals $\int_0^\infty \phi(w,t) \mathrm d w$ 
are convergent. On the other hand, step 1 of the algorithm outlined above requires to
numerically approximate this integral, and to this end, classical results from approximation theory 
\cite{Lubinsky} imply that such an algebraic decay does not admit a very fast convergence of
such numerical methods. Indeed, as the constant $q^{\mathrm C}$ is slightly larger than $q^{YA}$ for 
$\alpha \ge 1/2$ and (significantly) smaller for $\alpha < 1/2$, one may state that overall 
Chatterjee's method has more preferable properties from this point of view (although its 
properties are still far from good enough).
To the best of the author's knowledge, this is a feature shared by very many 
algorithms based on this type of approach. Therefore, one needs a relatively large number $K$ of quadrature nodes
in eq.\ \eqref{eq:qf1} to obtain an approximation with an acceptable accuracy (with the approaches known
so far, a common choice for $K$ is in the range between 200 and 500, cf.\ \cite{Ba2019,HSL2019}).
This number $K$ clearly has a strong influence on the constants implied in the $O$-term for the
computational complexity estimate. The main goal of this paper thus is to develop a 
method that is also based on the same fundamental idea but that leads to a function $\phi(w,t)$
which exhibits an exponential decay for large $w$. This behaviour is much more pleasant from
an approximation theoretic point of view because it allows to use well understood and rapidly 
convergent classical techniques like Gauss-Laguerre quadrature formulas. The hope behind this
idea is that the improved convergence behaviour will admit to use quadrature formulas as in
eq.\ \eqref{eq:qf1} with a significantly smaller number $K$ of nodes, so that the resulting algorithms 
can produce results with a comparable accuracy as the known methods in a much shorter amount of time
(that is still proportional to $N$ but with a significantly smaller implied constant).

\section{The New Diffusive Representation and its Properties}

Our main idea is based on the following result.

\begin{theorem}
	\label{thm:phinew}
	For given values $a \in \mathbb R$, $T > 0$ and $\alpha \in \mathbb R_+ \setminus \mathbb N$
	and a given function $y \in C^{\lceil \alpha \rceil}[a, a+T]$, let 
	\begin{equation}
		\label{eq:def-q}
		q_{\mathrm D} = \alpha - \lceil \alpha \rceil +1
	\end{equation}
	and
	\begin{equation}
		\label{eq:def-phinew}
		\phi^{\mathrm D}(w, t) 
			= (-1)^{\lfloor \alpha \rfloor} \frac{\sin \alpha \pi}{\pi} \mathrm e^{wq_{\mathrm D}}
				\int_a^t y^{(\lceil \alpha \rceil)}(\tau) \exp\left(-(t-\tau) \mathrm e^w\right) \mathrm d\tau
	\end{equation}
	for all $w \in \mathbb R$ and $t \in [a, a+T]$.
	Then, we have the following properties:
	\begin{enumerate}
	\item[(a)] The value $q_{\mathrm D}$ satisfies $0 < q_{\mathrm D} < 1$.
	\item[(b)] For any $w \in \mathbb R$, the function $\phi^{\mathrm D}(w, \cdot)$ solves the initial value
		problem 
		\begin{equation}
			\label{eq:ivp-phinew}
			\frac{\partial \phi^{\mathrm D}}{ \partial t} (w, t)
			= - \mathrm e^w \phi^{\mathrm D}(w, t) 
				+  (-1)^{\lfloor \alpha \rfloor} \frac{\sin \alpha \pi}{\pi} 
					\mathrm e^{wq_{\mathrm D}} y^{(\lceil \alpha \rceil)}(t),
			\quad
			\phi^{\mathrm D}(w, a) = 0
		\end{equation}
		for $t \in [a, a+T]$.
	\item[(c)] For any $t \in [a, a+T]$, 
		\begin{equation}
			\label{eq:phinew-caputo}
			D_a^\alpha y(t) = \int_{-\infty}^\infty \phi^{\mathrm D}(w, t) \mathrm d w.
		\end{equation}
	\item[(d)] For any $t \in [a, a+T]$, we have $\phi^{\mathrm D}(\cdot, t) \in C^\infty(\mathbb R)$.
	\item[(e)] For any $t \in [a, a+T]$, 
		\begin{equation}
			\label{eq:asymp-phi+}
			\phi^{\mathrm D}(w, t) = O(\mathrm e^{w (q_{\mathrm D}-1)}) 
			\quad \mbox{ as } w \to \infty
		\end{equation}
		and
		\begin{equation}
			\label{eq:asymp-phi-}
			\phi^{\mathrm D}(w, t) = O(\mathrm e^{w q_{\mathrm D}}) 
			\quad \mbox{ as } w \to -\infty.
		\end{equation}		
	\end{enumerate}
\end{theorem}

So, part (b) of Theorem \ref{thm:phinew} asserts that our function $\phi^{\mathrm D}$ solves an 
initial value problem of the same type as the previously considered functions, cf.\ \eqref{eq:ya-ode} or
\eqref{eq:sc-ode}. Moreover, according to part (c), by integrating this function with respect to $w$
we obtain the fractional derivative of the given function $y$, which is in analogy with the corresponding
equation \eqref{eq:diff-rep} for the known approaches mentioned above. Note that there is a marginal
difference between eqs.\ \eqref{eq:diff-rep} and \eqref{eq:phinew-caputo} in the sense that the latter
involves an integration over the entire real line whereas the former requires to integrate over the
positive half line only, but from the point of view of approximation (or quadrature) theory this does not 
introduce any substantial problems. (The index $\mathrm D$ in $\phi^{\mathrm D}$ and $q_{\mathrm D}$
can be interpreted to stand for ``doubly infinite integration range''.) 
Thus, in these respects the new model behaves in very much the 
same way as the known ones. The significant difference between the known approach and the new
one is evident from part (e) of the Theorem: It asserts (in view of the property of $q_{\mathrm D}$
shown in part (a)) that the integrand exhibits the desired exponential decay as $w \to \pm \infty$, 
thus allowing, in combination with the smoothness result of part (d), a much more efficient numerical integration.

\begin{proof}
	Part (a) is an immediate consequence of the definition of $q_{\mathrm D}$ given in eq.~\eqref{eq:def-q}.
	
	For part (b), we first note that the integrand in eq.\ \eqref{eq:def-phinew} is continuous by 
	assumption. Hence, the integral is zero for $t=a$ which implies that the initial condition given in
	eq.\ \eqref{eq:ivp-phinew} is correct. Also, a standard differentiation of the integral in the definition
	\eqref{eq:def-phinew} with respect to the parameter $t$ yields the differential equation.
	
	To prove (c), we recall from \cite[Proof of Theorem 3.18]{Di2010} that
	\[
		D_a^\alpha y(t) 
			= (-1)^{\lfloor \alpha \rfloor} \frac{\sin \alpha \pi} \pi 
				\int_a^t \int_0^\infty \frac{\mathrm e^{-z}} z \left( \frac z {x - \tau} \right)^{q_{\mathrm D}}
								y^{(\lceil \alpha \rceil)}(\tau) \mathrm d z \, \mathrm d \tau.
	\]
	The substitution $z = (x-\tau) \mathrm e^w$, combined with an interchange of the two integrations
	(that is admissible in view of Fubini's theorem), then leads to the desired result.
	
	Statement (d) directly follows from the definition \eqref{eq:def-phinew} of the function $\phi^{\mathrm D}$.
	
	Finally, we show that the estimates of (e) are true. To this end, let us first discuss what happens for
	$w \to + \infty$. Here, we can see that
	\[
		\phi^{\mathrm D}(w, t) =(-1)^{\lfloor \alpha \rfloor} \frac{\sin \alpha \pi}{\pi} ( I_1 + I_2 )
	\]
	where
	\begin{align*}
		| I_1 |
		&= \left| \mathrm e^{w q_{\mathrm D}} 
			\int_{t - w \exp(-w)}^t  y^{(\lceil \alpha \rceil)}(\tau) 
						\exp\left(-(t-\tau) \mathrm e^w\right) \mathrm d\tau \right| \\
		& \le \| y^{(\lceil \alpha \rceil)} \|_{L_\infty[a, a+T]} \mathrm e^{w q_{\mathrm D}} 
				\left| \int_{t - w \exp(-w)}^t 
						\exp\left(-(t-\tau) \mathrm e^w\right) \mathrm d\tau \right| \\
		& \le \| y^{(\lceil \alpha \rceil)} \|_{L_\infty[a, a+T]} \mathrm e^{w q_{\mathrm D}} 
				\mathrm e^{-w} \left[ 1 - \mathrm e^{-w}  \right] 
		 < \| y^{(\lceil \alpha \rceil)} \|_{L_\infty[a, a+T]} \mathrm e^{w (q_{\mathrm D} - 1)} 
	\end{align*}
	and
	\begin{align*}
		| I_2 |
		&= \left| \mathrm e^{w q_{\mathrm D}} 
			\int_a^{t - w \exp(-w)}  y^{(\lceil \alpha \rceil)}(\tau) 
						\exp\left(-(t-\tau) \mathrm e^w\right) \mathrm d\tau \right| \\
		&\le \mathrm e^{w q_{\mathrm D}} 
			\max_{\tau \in [a, t - w \exp(-w)]} \exp\left(-(t-\tau) \mathrm e^w\right) 
			\int_a^{t - w \exp(-w)} | y^{(\lceil \alpha \rceil)}(\tau) |  \mathrm d\tau \\
		&\le \mathrm e^{w q_{\mathrm D}} 
			\mathrm e^{-w}
			\int_a^{a+T} | y^{(\lceil \alpha \rceil)}(\tau) |  \mathrm d\tau 
		 = \mathrm e^{w (q_{\mathrm D}-1)} 
			\int_a^{a+T} | y^{(\lceil \alpha \rceil)}(\tau) |  \mathrm d\tau 
	\end{align*}
	which shows the desired result \eqref{eq:asymp-phi+} in this case; in particular,
	the upper bound decays exponentially for $w \to \infty$ because $q_{\mathrm D}<1$. 
	Regarding the behaviour for $w \to -\infty$, we start from the representation \eqref{eq:def-phinew} 
	and apply a partial integration. 
	This yields, taking into consideration that $t \ge a$, that
	\begin{align*}
		| \phi^{\mathrm D}(w, t) |
		&= \frac{|\sin \alpha \pi|}{\pi} \mathrm e^{wq_{\mathrm D}}
			 \Big |  \exp \left(-(x-\tau) \mathrm e^w\right)  
			 		y^{(\lceil \alpha \rceil - 1)}(\tau)  {\big|}_{\tau=a}^{\tau=t} \\
			 		& \qquad \qquad \qquad \qquad
			 		- \mathrm e^w \int_a^t  \exp \left( -(t-\tau) \mathrm e^w \right) 
			 			y^{(\lceil \alpha \rceil - 1)}(\tau) \mathrm d \tau \Big | \\
		&\le  \frac{|\sin \alpha \pi|}{\pi} \mathrm e^{wq_{\mathrm D}} 
			\left| 
				y^{(\lceil \alpha \rceil - 1)}(t) 
				- y^{(\lceil \alpha \rceil - 1)}(a) \exp\left(-(t-a) \mathrm e^w \right)
			\right| \\
		& \phantom{\le} \quad {} + \frac{|\sin \alpha \pi|}{\pi} \mathrm e^{w q_{\mathrm D}}
			\| y^{(\lceil \alpha \rceil - 1)} \|_{L_\infty[a, a+T]} 
				\left| \mathrm e^w \int_a^t \exp \left( -(t-\tau) \mathrm e^w \right) \mathrm d \tau \right| \\
		&\le  \frac{|\sin \alpha \pi|}{\pi} 
				\| y^{(\lceil \alpha \rceil - 1)} \|_{L_\infty[a, a+T]} \mathrm e^{wq_{\mathrm D}} 
				\left( 2 + 1 - \exp \left( -(t-a) \mathrm e^w \right) \right) \\
		&\le  3 \frac{|\sin \alpha \pi|}{\pi} 
				\| y^{(\lceil \alpha \rceil - 1)} \|_{L_\infty[a, a+T]} \mathrm e^{wq_{\mathrm D}} ,				
	\end{align*}
	thus proving the relation \eqref{eq:asymp-phi-} and demonstrating, in view of $q_{\mathrm D} > 0$,
	that $\phi^{\mathrm D}(w, t)$ decays to zero exponentially as $w \to -\infty$.
\end{proof}

\section{The Complete Numerical Method}

Based on Theorem \ref{thm:phinew}---in particular, using the properties shown in parts (d) and (e)---we 
thus proceed as follows to obtain the required approximation of $D_a^\alpha y(t_j)$, $j = 1, 2, \ldots, N$.
Splitting up the integral from eq.\ \eqref{eq:phinew-caputo} into the integrals over the negative 
and over the positive half line, respectively, and introducing some obvious substitutions, we notice that
\begin{align*}
	\int_{-\infty}^\infty \phi^{\mathrm D}(w, t) \mathrm d w
	&= \frac 1 {q_{\mathrm D}} \int_{0}^\infty \mathrm e^{-u} \mathrm e^{u} 
										\phi^{\mathrm D}(-u/q_{\mathrm D}, t) \mathrm d u \\
	& \phantom{=} \quad {}
		+ \frac 1 {1-q_{\mathrm D}} \int_{0}^\infty \mathrm e^{-u} \mathrm e^{u} 
										\phi^{\mathrm D}(u/(1-q_{\mathrm D}), t) \mathrm d u.
\end{align*}
Therefore, using 
\begin{equation}
	\label{eq:def-phitilde}
	\hat \phi^{\mathrm D}(u, t) 
		:=  \mathrm e^u \left( \frac 1 {q_{\mathrm D}} \phi^{\mathrm D}(-u/q_{\mathrm D}, t)
						+  \frac 1 {1-q_{\mathrm D}} \phi^{\mathrm D}(u/(1-q_{\mathrm D}), t) \right),
\end{equation}
we find that
\begin{equation}
	D^\alpha_a y(t) 
	= \int_{-\infty}^\infty \phi^{\mathrm D}(w, t) \mathrm d w 
	= \int_0^\infty \mathrm e^{-u} \hat \phi^{\mathrm D}(u, t) \mathrm d u 
	 \approx Q^{\mathrm{GLa}}_K [  \hat \phi^{\mathrm D}(\cdot, t)]
\end{equation}
where
\[
	 Q^{\mathrm{GLa}}_K [f] = \sum_{k=1}^K a^{\mathrm{GLa}}_k f(x^{\mathrm{GLa}}_k)
\]
is the $K$-point Gauss-Laguerre quadrature formula, i.e.\ the Gaussian quadrature formula for
the weight function $\mathrm e^{-u}$ on the interval $[0, \infty)$ \cite[Sections 3.6 and 3.7]{DR}.
For the sake of simplicity, we have chosen to omit from our notation for the nodes $x^{\mathrm{GLa}}_k$
and the weights $a^{\mathrm{GLa}}_k$ of the Gauss-Laguerre quadrature formula the fact that
these quantities depend on the total number $K$ of quadrature nodes.
From \cite [p.\ 227]{DR} and our Theorem \ref{thm:phinew} above, we can immediately
conclude the following result:

\begin{theorem}
	\label{thm:conv-gl}
	Under the assumptions of Theorem \ref{thm:phinew}, we have
	\[
		\lim_{K \to \infty} Q^{\mathrm{GLa}}_K [ \hat \phi^{\mathrm D}(\cdot, t) ] 
		= D^\alpha_a y(t) 
	\]
	for all $t \in [a, a+T]$.
\end{theorem}

For a given number $K$ of quadrature points, it is known that the 
nodes $x^{\mathrm{GLa}}_k$, $k = 1, 2, \ldots, K$, are the zeros of the 
Laguerre polynomial $L_K$ of order $K$, and the associated weights are given by
\[
	 a^{\mathrm{GLa}}_k =  \frac{x^{\mathrm{GLa}}_k}{[L_{K+1}(x_k)]^2},
\]
cf., e.g., \cite[p.\ 223]{DR}. (In our definition of the Laguerre polynomials, the normalization is such that 
$\int_0^\infty \mathrm e ^{-x} (L_K(x))^2 \mathrm d x = 1$.) 
From \cite[eqs.\ (6.31.7), (6.31.11) and (6.31.12)]{Sz} we know that,
at least for $K \ge 3$,
\[
	\frac{2.89}{2K+1} < x^{\mathrm{GLa}}_1 < \frac 3 {2 K}
	\quad \mbox{ and } \quad
	2 K < x^{\mathrm{GLa}}_K < 4K + 3.
\]

We are now in a position to describe the method for the numerical computation 
of $D^\alpha_a y(t_j)$, $j = 1, 2, \ldots, N$, that we propose. In this algorithm, 
the symbol $\phi_k$ is used to denote the approximate value of 
$\phi^{\mathrm D}(x^{\mathrm{GLa}}_k, t_j)$ for the current time step. i.e.\ for the 
currently considered value of $j$.
Steps 1 and 2 here are merely preparatory in nature; the core of the algorithm is step 3.

\begin{quote}
	Given the initial point $a$, the order $\alpha$, the grid points $t_j$, $j = 1, 2, \ldots, N$ 
	and the number $K \in \mathbb N$ of quadrature nodes,
	\begin{enumerate}
	\item Set $q_{\mathrm D} \mapsfrom \alpha - \lceil \alpha \rceil + 1$.
	\item For $k = 1, 2, \ldots, K$:
		\begin{enumerate}
		\item compute the Gauss-Laguerre nodes $x^{\mathrm{GLa}}_k$ 
			and the associated weights $a^{\mathrm{GLa}}_k$,
		\item define the auxiliary quantities $w_k \mapsfrom - x^{\mathrm{GLa}}_k / q_{\mathrm D}$
			and $\tilde w_k \mapsfrom x^{\mathrm{GLa}}_k / (1 - q_{\mathrm D})$,
		\item set $\phi_k \mapsfrom 0$ and $\tilde \phi_k \mapsfrom 0$
			(to represent the initial condition of the differential equation \eqref{eq:ivp-phinew} 
			for $t = t_0 = a$).
		\end{enumerate}
	\item For $j = 1, 2, \ldots, N$:
		\begin{enumerate}
		\item Set $h \mapsfrom t_j - t_{j-1}$.
		\item For $k = 1, 2, \ldots, K$: 
			\begin{subequations}
			\begin{enumerate}
			\item update the value $\phi_k$ by means of solving the associated differential equation 
				\eqref{eq:ivp-phinew} with, e.g., the backward Euler method, viz.\ 
				\begin{equation}
					\label{eq:bweuler}
					\phi_k \mapsfrom \frac 1 {1 + h \mathrm e^{w_k}} 
										\left( \phi_k + h (-1)^{\lfloor \alpha \rfloor} 
												\frac{\sin \alpha \pi}{\pi} 
												\mathrm e^{w_k q_{\mathrm D}} 
												y^{(\lceil \alpha \rceil)}(t_{j})
										\right)
				\end{equation}
				(note that the index $k$ used here is not the time index);
			\item similarly, update the value $\tilde \phi_k$ by 
				\begin{equation}
					\label{eq:bweuler-tilde}
					\tilde \phi_k \mapsfrom \frac 1 {1 + h \mathrm e^{\tilde w_k}} 
										\left( \tilde \phi_k + h (-1)^{\lfloor \alpha \rfloor} 
												\frac{\sin \alpha \pi}{\pi} 
												\mathrm e^{\tilde w_k q_{\mathrm D}} 
												y^{(\lceil \alpha \rceil)}(t_{j})
										\right).
				\end{equation}
			\end{enumerate}
			\end{subequations}
		\item Compute the desired approximate value for $D^\alpha_a y(t_j)$ using the formula
			\[
				D^\alpha_a y(t_j) 
					= \sum_{k=1}^K a^{\mathrm{GLa}}_k 
									\exp(x^{\mathrm{GLa}}_k) 
									\left( \frac 1 {q_{\mathrm D}} \phi_k
										+ \frac 1 {1 - q_{\mathrm D}} \tilde \phi_k
									\right) .
			\]
		\end{enumerate}
	\end{enumerate}
\end{quote}

The main goal of this paper is to develop a diffusive representation that can be numerically handled in 
a more efficient way than traditional formuals. Therefore, our work concentrates on the aspects 
related to the integral, i.e.\ on the properties of the integrand and on the associated numerical 
quadrature. The solution of the 
differential equation is not in the focus of our work; we only use some very simple (but nevertheless
reasonable) methods here. Our specific choice is based on the observation that the magnitude of the 
constant factor with which the unkonwn function $\phi(w, \cdot)$ on the right-hand side of 
\eqref{eq:ivp-phinew} is multiplied is such that an A-stable method should be used \cite{HW}.
Therefore, as the simplest possible choice among these methods, we have suggested the 
backward Euler method in our description given above.
Alternatively, one could, e.g., use the trapezoidal method which is also A-stable. This would mean that
the formulas given in eqs.~\eqref{eq:bweuler} and \eqref{eq:bweuler-tilde} would have to be
replaced by
\begin{subequations}
\begin{align}
	\label{eq:trap}
		\phi_k &\mapsfrom \frac 1 {1 + h \mathrm e^{w_k} /2} 
							\Bigg(  \left(1 - \frac h 2 \mathrm e^{w_k} \right) \phi_k \\
			& \phantom{mapsfrom} {} \qquad \qquad
								+ \frac h 2  (-1)^{\lfloor \alpha \rfloor} 
									\frac{\sin \alpha \pi}{\pi} 
									\mathrm e^{w_k q_{\mathrm D}} 
									(y^{(\lceil \alpha \rceil)}(t_{j}) + y^{(\lceil \alpha \rceil)}(t_{j-1}))
							\Bigg)
							\nonumber
\end{align}
and
\begin{align}
	\label{eq:traptilde}
		\tilde \phi_k &\mapsfrom \frac 1 {1 + h \mathrm e^{\tilde w_k} /2} 
							\Bigg(  \left(1 - \frac h 2 \mathrm e^{\tilde w_k} \right) \tilde \phi_k \\
			& \phantom{mapsfrom} {} \qquad \qquad
								+ \frac h 2  (-1)^{\lfloor \alpha \rfloor} 
									\frac{\sin \alpha \pi}{\pi} 
									\mathrm e^{\tilde w_k q_{\mathrm D}} 
									(y^{(\lceil \alpha \rceil)}(t_{j}) + y^{(\lceil \alpha \rceil)}(t_{j-1}))
							\Bigg)
							\nonumber
\end{align}
\end{subequations}
respectively. In the following section, we shall report the results of our numerical experiments
for both variants.

\begin{remark}
	From a formal point of view, eqs. \eqref{eq:bweuler} and \eqref{eq:bweuler-tilde} have
	exactly the same structure. From a numerical perspective, however, there is a significant
	difference between them that needs to be taken into account when implementing the
	algorithm in finite-precision arithmetic: In view of the definitions of the quantities $w_k$ and
	$\tilde w_k$ given in step 2b of the algorithm and the facts that the Gauss-Laguerre nodes
	$x_k^{\mathrm{GLa}}$ are strictly positive for all $k$ and that $q_{\mathrm D} \in (0,1)$,
	it is clear that $w_k < 0$ for all $k$, and hence the powers $\mathrm e^{w_k}$ and 
	$\mathrm e^{w_k q_{\mathrm D}}$ that occur in eq.\ \eqref{eq:bweuler} are always in the
	interval $(0,1)$. It may be, if $|w_k|$ is very large, that the calculation of $\mathrm e^{w_k}$ 
	in IEEE arithmetic results in an underflow, but this number can then safely be replaced by $0$ 
	without causing any problems. Therefore, eq.\ \eqref{eq:bweuler} can be implemented directly
	in its given form. On the other hand, using the same arguments we can see that $\tilde w_k > 0$
	for all $k$, and indeed (at least if $k$ is large and/or $q_{\mathrm D}$ is close to $1$) 
	$\tilde w_k$ may be so large that the computation of $\mathrm e^{\tilde w_k}$ results in a fatal 
	overflow. For this reason, in a practical implementation, eq.\ \eqref{eq:bweuler-tilde} should
	not be used in its form given above but in the equivalent form
	\begin{equation}
		\label{eq:bweuler-tilde2}
		\tag{14c}
		\tilde \phi_k \mapsfrom \frac {\mathrm e^{-\tilde w_k}} {\mathrm e^{-\tilde w_k} + h} \tilde \phi_k 
						+ h (-1)^{\lfloor \alpha \rfloor} \frac{\sin \alpha \pi}{\pi}
						  \frac {\mathrm e^{\tilde w_k (q_{\mathrm D} - 1)}} {\mathrm e^{-\tilde w_k} + h} 
							y^{(\lceil \alpha \rceil)}(t_{j})
	\end{equation}
	that avoids all potential overflows.
	
	Evidently, an analog comment applies to eqs.\ \eqref{eq:trap} and \eqref{eq:traptilde}.
\end{remark}

\section{Experimental Results and Conclusion}

In \cite{Di2021}, we have reported some numerical results illustrating the convergence behaviour
of the RISS method proposed by Hinze et al.\ \cite{HSL2019}. Here now we 
present similar numerical results obtained with the new algorithm. 
A comparison with the corresponding data shown in \cite{Di2021} reveals that, in many cases, 
our new method requires a smaller number of quadrature nodes than the RISS approach
(with otherwise identical parameters) to obtain approximations of a similar quality. 

\begin{figure}
	\vskip-0.4cm
	\centering
	\includegraphics[width=0.7\textwidth]{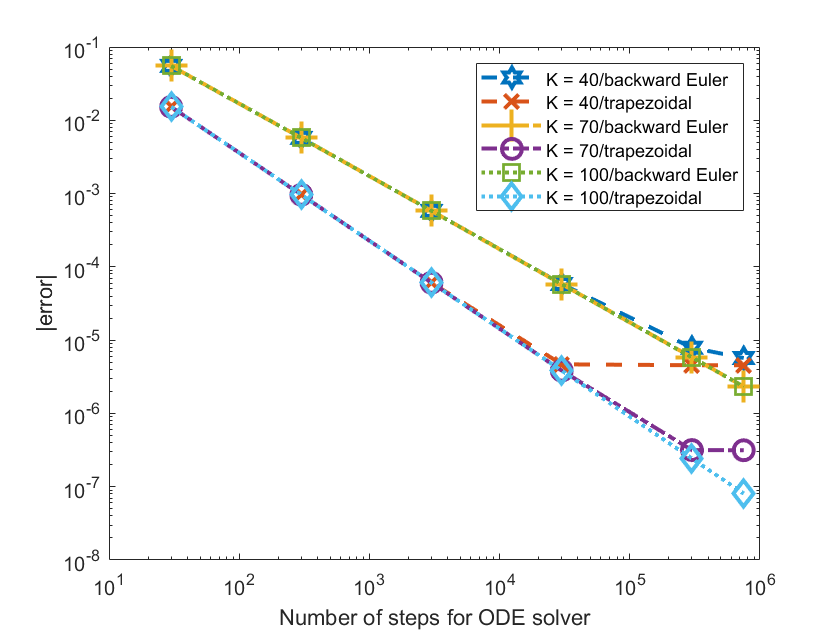}
	\caption{\label{fig:ex1}Maximal errors for the calculation of $D_0^{0.4} y(t)$ with $t \in [0, 3]$ for
		$y(t) = t^{1.6}$ using different step sizes for the ODE solver and different numbers of 
		quadrature nodes.}
\end{figure}

A typical result is shown in Figure \ref{fig:ex1} where we have numerically computed 
the Caputo derivative of order $0.4$ of the function $y(t) = t^{1.6}$ over the interval $[0,3]$.
The calculations have been performed on an equispaced grid for the interval $[0,3]$ 
with various different step sizes (i.e.\ with different numbers of grid points) and different 
choices of the number $K$ of quadrature nodes. Both the backward Euler and the trapezoidal scheme
have been tried as the ODE solvers. The figure exhibits the maximal absolute error over all grid 
points.

The findings of this example can be summarized as follows:
\begin{itemize}
\item The trapezoidal method clearly leads to a more accurate approximation than the 
	backward Euler method. Obviously, in view of the trapezoidal method's higher
	convergence order, this behaviour is exactly what would have been expected.
\item The number of quadrature points, i.e.\ our parameter $K$, only has a very small influence
	on the overall error. Therefore, one can afford to work with a relatively small value of $K$,
	thus significantly reducing the computational cost, without a substantial loss of accuracy.
\item A comparison of the results for the trapezoidal method shows that a certain kind of saturation 
	is reached at an error level of $4.5\cdot10^{-6}$ for $K=40$, i.e.\ we do not achieve 
	a better accuracy even if we continue to decrease the step size for the ODE solver.
	This is an indication that this level reflects the contribution of the total error caused
	by the quadrature formula. If a smaller error is required, one therefore needs to use more 
	quadrature nodes. For example, choosing $K = 70$ leads to a saturation level
	of approximately $3.2\cdot10^{-7}$. This indicates that the saturation level
	might be proportional to $K^{-0.6}$, leading to the conjecture that the exponent 
	of $K$ in this expression could be related to the smoothness properties of the 
	function $y$ (note that the function $y'$ that appears in the formulas which describe our
	algorithm satisfies a Lipschitz condition of order $0.6$).
	
	The fact that this phenomenon is hardly visible if the backward Euler method is used
	is due to the fact that this ODE solver has a larger error which only just about reaches this
	range for the chosen step sizes. It would be possible to more clearly observe a similar behaviour
	if the step sizes were reduced even more.
\end{itemize}
We have also used a number of other test cases; the behaviour has usually been very similar. 
Also, the findings of \cite{Di2021} for a significantly different method based on a related fundamental
approach point into the same direction. In our future work, we will attempt to provide a thorough
analysis of the approximation properties of methods of this type that should confirm the 
experimental results.

\end{document}